\documentstyle[12pt]{article}
\newcommand{\ben}{\begin{enumerate}}
\newcommand{\een}{\end{enumerate}}
\newcommand{\be}{\begin{equation}}
\newcommand{\ee}{\end{equation}}
\newcommand{\bea}{\begin{eqnarray}}
\newcommand{\eea}{\end{eqnarray}}
\newcommand{\bc}{\begin{center}}
\newcommand{\ec}{\end{center}}

\newcommand{\la}{\leftarrow}

\newcommand{\ba}{\begin{array}}
\newcommand{\ea}{\end{array}}

\begin{document}
\begin{center}
{\normalsize\bf
The consistent reduction of the differential calculus on the
quantum group $GL_{q}(2,C)$ to the differential calculi on its
subgroups and $\sigma$-models on the quantum group manifolds
$SL_{q}(2,R)$, $SL_{q}(2,R)/U_{h}(1)$, $C{q}(2|0)$ and infinitesimal
transformations}\\
\vspace{0.20cm} V.D.  Gershun \vspace{0.20cm}\\
{\it  Institute of Theoretical Physics, NSC Kharkov Institute of
Physics and Technology, }\\ {\it 310108 Kharkov, Ukraine }\\
\end{center}
\begin{abstract}

Explicit construction of the second order left differential
calculi on the quantum group and its subgroups are obtained with the
property of the natural reduction: the differential calculus on the
quantum group  $GL_q(2,C)$ has to contain the 3-dimensional differential
calculi on the quantum subgroup $SL_q(2,C)$, the differential calculi
on the Borel subgroups $B_{L}^{(2)}(C)$, $B_{U}^{(2)}(C)$ of the lower and
of the upper triangular matrices, on the quantum subgroups $U_{q}(2)$,
$SU_{q}(2)$, $Sp_{q}(2,C)$, $Sp_{q}(2)$, $T_{q}(2,C)$, $B_{L}(C)$,
$B_{U}(C)$, $U_{q}(1)$, $Z_{-}^{(2)}(C)$, $Z_{+}^{(2)}(C)$ and on the their
real forms. The classical limit ($q\to 1$)
of the  left  differential calculus is the nondeformed
differential calculus. The differential calculi on the Borel subgroups
$B_{L}(C)$, $B_{U}(C)$ of the $SL_{q}(2,C)$ coincide with two solutions of
Wess-Zumino differential calculus on the quantum plane $C_q(2|0)$.

The spontaneous breaking symmetry in the WZNW model
with $SL_{q}(2,R)$ quantum group symmetry over two-dimensional
nondeformed Minkovski space and in the $\sigma$-models with
${SL_{q}(2,R)/U_{p}(1)}$, $C_{q}(2|0)$ quantum group symmetry is
considered. The Lagrangian formalism over the quantum group manifolds is
discussed. The variational calculus on the $SL_{q}(2,R)$ group manifold
is obtained. The classical solution of $C_{q}(2|0)$ {$\sigma$}-model is
obtained.  \end{abstract}

\noindent{\bf I. Introduction}

The left-(right-) invariant differential calculus on the $GL_{q}(2,C)$,
matched with differential calculi on the $SL_{q}(2,C)$ and on the
$C_{q}(2|0)$, was constructed in [1,2], the differential calculus on the
$GL_{q}(2,C)$, matched with the differential calculi on its subgroups, was
constructed in [3]. In this paper we shortly remember the basic states of
the left-invariant consistent differential calculus and on this
base we consider the two-dimensional WZNW model on the $SL{q}(2,R)$ group
over nondeformed $M^{(1,1)}$ space. We obtained the expression for the
background gravity and antisymmetric fields in the $\sigma$-model approach
to a relativistic string. Also, we considered the spontaneous breaking
symmetry in this model and in the $\sigma$-models on the quantum group
manifolds $SL_{q}(2,R)/U_{h}(1)$ and on the $C_{q}(2|0)$. We considered
the infinitesimal transformations on the $SL_{q}(2,R)$ and obtained the
commutation relations of the variational calculus on this group.

\noindent{\bf II. Differential calculus on the $GL_{q}(2,C)$ group.}

The matrix quantum group [1] $G = GL_{q}(2,C)$ is defined by the
q-commutation relations (C.R.) of its group parameters. Let $g \in
GL_{q}(2,C)$
\be
g = \left(\begin{array}{cc}
g_{1}^{1} & g_{2}^{1} \\ g_{1}^{2} & g_{2}^{2}
\end{array}\right )  =
\left(\begin{array}{cc} a & b \\ c & d
\end{array} \right ),\;\;\;
R_{kl}^{ij}g_m^kg_n^l = g_l^jg_k^iR_{mn}^{kl}
\ee
The $Rgg$-relations in components have the standard form
\be
\begin{array}{llll} ab = qba, & bd = qdb, & bc = cb, &
\lambda=q-\frac{1}{q} \\ ac = qca, & cd = qdc, & ad = da + \lambda bc
\end{array} \ee The quantum determinant $D_{q}$ is the central
element and commutes with the group elements \be
det_{q}g = D_{q} = ad - qbc = da - \frac{1}{q}cb,\;\;\; D_{q}g_{i}^{k} =
g_{i}^{k}D_{q} \ee
The left-side exterior derivative of any function $f$ of group parameters
is \be df = \left( f{\frac{\stackrel{\la }{\partial }}
{\partial
g_i^k}} \right ) dg_i^k \ee The differential calculus on a quantum
group is
defined by C.R. between the group parameters and theirs differentials, or
between the group parameters and Maurer-Cartan 1-forms \be
\omega =
g^{-1}dg = \left(\begin{array}{cc}\omega_1^1 & \omega_2^1\\ \omega_1^2 &
\omega_2^2\end{array}\right ) = \left(\begin{array}{cc}\omega ^1 & \omega
^2\\ \omega ^3 & \omega ^4 \end{array}\right ) \ee \noindent
The basic
states of the left-invariant differential calculus are following:

{\tt 1.}Following [5], we look for C.R. between 1-forms and
group parameters in the form
\be
\omega^{\alpha}_{\beta}T^i_k = T^i_m
\left( F^m_k\right ) ^{\alpha\gamma}_{\delta\beta}
\omega^{\delta}_{\gamma},\;\;\; F^m_k = \left(\begin{array}{cc} A & B\\ C
& D \end{array}\right ) \ee and $A, B, C, D$ are $4\times 4$
unknown matrices with complex entries.

{\tt 2.} The requirement
\be [\omega^i_j, \left (R^{\alpha n}_{\beta m} T^{\beta}_
{\delta} T^m_k -
T^n_m T^{\alpha}_{\beta} R^{\beta m}_{\delta k}\right )] = 0 \ee
leads to
the equations for the matrices $A,B,C,D$ \be R^{\alpha n}_{\beta m}
F^{\beta}_{\delta} F^m_k = F^n_m F^{\alpha}_{\beta} R^{\beta m}_{\delta k}
\ee $A,B,C,D$ are the representation of $a,b,c,d.$

{\tt 3.} For
the quantum determinant $D_{q}$ we define the quantum trace by the
definition \be dD_{q} = D_{q}{\mbox{Tr}_q}\omega \ee
We obtained \be
\begin{array}{lll} dD_{q} = D_{q}\left(\omega^1 + A^4_l\omega^l -
\frac{1}{q}B^3_l\omega^l\right ) & = D_{q}\left(\omega^4+D^1_l\omega^l -
qC^2_l\omega^l\right )\\ = D_{q}\left(\omega^1+A^4_l\omega^l -
\frac{1}{q}C^2_l\omega^l\right ) & = D_{q}\left(\omega^4+D^1_l\omega^l -
qC^2_l\omega^l\right ) \end{array} \ee and additional conditions
\be
\omega^2 + B^1_l\omega^l - qA^2_l\omega^l = 0 , \;\;\;
\omega^3 + C^4_l\omega^l - \frac{1}{q} D^3_l\omega^l = 0
\ee

{\tt 4.} The C.R. between 1-forms and the quantum determinant
are
\be
\omega^i_kD_{q} = D_{q}\left(AD-qBC\right )^i_m\omega^m_k
\ee

{\tt 5.} Differentiating the commutation relations between 1-forms
and group parameters and using the Maurer--Cartan equation $d\omega = -
\omega\wedge\omega$ we have the $q$-deformation of the algebra of
the $\omega$-forms:
\be
\left(F^f_m\right )^{\alpha\varphi}_{\tau\gamma}
\left(F^m_k\right )^{\gamma\rho}_{\sigma\beta}
\omega^{\tau}_{\varphi} \omega^{\sigma}_{\rho} +
\left(F^f_{\sigma}\right )^{\alpha\varphi}_{\tau\beta}
\omega^{\tau}_{\varphi} \omega^{\sigma}_{k} +
\left(F^{\varphi}_{k}\right )^{\alpha\rho}_{\sigma\beta}
\omega^{f}_{\varphi} \omega^{\sigma}_{\rho} -
\left(F^{f}_{k}\right )^{\alpha\rho}_{\tau\beta}
\omega^{\tau}_{\sigma} \omega^{\sigma}_{\rho}
= 0
\ee
\noindent All of this conditions are common for any differential calculus.

{\tt 6.}The last condition is the basic condition,
which differs the matched differential calculus from the bicovariant
differential calculus [6]. Instead of bicovariance condition we require
\be
\omega^kD_{q} = D_{q}\omega^k
\ee
which leads to the equation $AD - qBC = 1$.
To solve the system of equations for the matrix $F$ satisfying to this
conditions we considered the following representation of algebra (8)
\footnote{ Note,that the solutions for the matrix $F$ with $B,C \neq 0$ ,
considered by author, lead to the differential calculi which have not
classical limit $(q\to 1)$.}
\be B = C = 0\;\;\;\;\;\;AD = 1.  \ee
We obtained the following expressions for the quantum trace
\be
dD_{q} = D_{q}\left[D^1_1 \omega^1 + (1+D^1_4)\omega^4\right] =
D_{q}\left[(1 + A^4_1)\omega^1 + A^4_4 \omega^4\right],
\ee
where
we used the composite indexes 1,2,3,4.
Finally, we found the matrices $A,D$ in terms of the independent variables
$\beta = D^1_1,\;\alpha = A^4_4$
\medskip
\be
A = \left(\begin{array}{cccc}
1-\alpha + \frac{\alpha}{\beta} & 0 & 0 & \frac{(1 -
\alpha)\alpha}{\beta}\\ 0 & \frac{1}{q} & 0 & 0 \\ 0 & 0 & \frac{1}{q} &
0\\ \beta-1 & 0 & 0 & \alpha\\ \end{array}\right ), \;\;\;
D = \left(\begin{array}{cccc}
\beta & 0 & 0 & \alpha - 1 \\
0 & q & 0 & 0 \\
0 & 0 & q & 0 \\
\frac{(1 - \beta)\beta}{\alpha} & 0 & 0 & 1 - \beta + \frac{\beta}{\alpha}
\end{array}\right )
\ee
We found, that the algebra of the  1-forms can be decomposed on the algebra
$SL_q(2,C)$ and the $U(1)$ subalgebra only for unique value of the
parameters
\be
\alpha = \frac{2}{1 + q^2}, \;\;\;\;\;\;
\beta = \frac{2q^2}{1 + q^2}\ee In this case the commutation
relations between the parameters and the 1-forms are diagonalized after
the choice of the new basis of 1-forms \be
\bar{\omega^1} = \frac{2}{q
+ 1/q}(q\omega^1 + \frac{1}{q}\omega^4) = {\mbox{Tr}_q}\omega, ~~~
\bar{\omega^4} = \frac{1}{1 + q^2} (\omega^1 - \omega^4)
\ee
\be
\begin{array}{lcl}
\bar{\omega^1}a = a\bar{\omega^1} & \bar{\omega^1}d = d\bar{\omega^1}\\
\bar{\omega^4}a = q^{-2}a\bar{\omega^4} &
\bar{\omega^4}d = q^2d\bar{\omega^4},\\ \omega^2a = q^{-1}a\omega^2 &
\omega^2d = qd\omega^2,\\ \omega^3a = q^{-1}a\omega^3 &
\omega^3d = qd\omega^3 \end{array} \ee The other relations are
found by making the interchanges $a\leftrightarrow c,d\leftrightarrow b$.
The commutation relations between 1-forms can be written as
\be
\begin{array}{ll}
{\bar{\omega^1}\omega^2 + \omega^2\bar{\omega^1} = 0} , & \;\;
{\bar{\omega^1}\omega^3 + \omega^3\bar{\omega^1} = 0} \cr
{q^2\omega^2\bar{\omega^4} + q^{-2}\bar{\omega^4}\omega^2 = 0} , & \;\;
{q^2\bar{\omega^4}\omega^3 + q^{-2}\omega^3\bar{\omega^4} = 0} \cr
{\bar{\omega^1}\bar{\omega^4} + \bar{\omega^4}\bar{\omega^1} = 0} , & \;\;
{\omega^2\omega^3 + q^{-2}\omega^3\omega^2 = 0}
\end{array}
\ee
\be
\begin{array}{ll}
{(\omega^1)^2 = (\omega^2)^2 = (\omega^3)^2 = (\omega^4)^2 = 0}
\end{array}
\ee
The algebra of the vector fields
${\stackrel{\la}{\nabla_{k}}}$ has the form
\be
\begin{array}{cc}
{\stackrel{\la}{\nabla_{3}}}{\stackrel{\la}{\nabla_{2}}} -
q^2{\stackrel{\la}{\nabla_{2}}}{\stackrel{\la}{\nabla_{3}}}
 = {\stackrel{\la}{\nabla_{1}}} &
q^2{\stackrel{\la}{\nabla_{2}}}{\stackrel{\la}{\nabla_{1}}}
- q^{-2}{\stackrel {\la}{\nabla_{1}}}{\stackrel{\la}
{\nabla_{2}}}
 = (1 + q^2){\stackrel{\la}{\nabla_{2}}} \cr
q^2{\stackrel{\la}{\nabla_{1}}}{\stackrel{\la}{\nabla_{3}}}
- q^{-2}{\stackrel {\la}{\nabla_{3}}}{\stackrel{\la}
{\nabla_{1}}}
 = (1+q^2){\stackrel{\la}{\nabla_{3}}} &
[{\stackrel{\la}{\nabla_{4}}},
{\stackrel{\la}{\nabla_{k}}}] = 0 \cr
\end{array}
\ee
After the mapping
\be {\stackrel{\la}{\nabla_{1}}},
{\stackrel{\la}{\nabla_{2}}},
{\stackrel{\la}{\nabla_{3}}},
{\stackrel{\la}{\nabla_{4}}} \to H, T_{2}, T_{3}, N \ee
\be
\begin{array}{cc}
{\stackrel{\la}{\nabla_{1}}} =
\frac{1 - q^{ - 2H}}{1 - q^{ - 2}}, &
\;\;\;{\stackrel{\la}{\nabla_{2}}} = q^{- H/2}T_{2},\cr
{\stackrel{\la}{\nabla_{4}}} = N ,& \;\;\;
{\stackrel{\la}{\nabla_{3}}} = q^{ - H/2}T_{3}
\end{array}
\ee
we obtain the algebra $ U_q GL(2,C) $ in the form of the Drinfeld--Jimbo
algebra.
\be\begin{array}{cc}
[T_{3}, T_{2}] = \frac{q^H - q^{ - H}}{q - q^{- 1}}, &
[H , T_{3}] = 2T_{3} \cr
[H , T_{2}] = - 2T_{2} & {[N , H] = [N , T_{2}] = [N , T_{3}] = 0}
\end{array}
\ee Using the commutation relations between the parameters and
1-forms we obtained the C.R. between the parameters and its left
differentials.
\be
\begin{array}{l}
da\;a = q^{- 2}\;a\;da +
\frac{q^2 - 1}{2q^2}\;a^2\;{\mbox{Tr}_q}\omega,
\cr
dc\;c = q^{-2}\;c\;dc +
\frac{q^2 - 1}{2q^2}\;c^2\;{\mbox{Tr}_q}\omega,\cr
da\;c = q^{-1}\;c\;da +
\frac{q^2 - 1}{2q^2}\;a\;c\;{\mbox{Tr}_q}\omega,
\cr
dc\;a=q^{-1}\;a\;dc + (q^{-2} - 1)\;c\;da + \frac{q^2 - 1}{2q^2}\;c\;a\;
{\mbox{Tr}_q}\omega,
\cr
db\;b = q^{2}\;b\;db +
\frac{1 - q^2}{2}\;b^2\;{\mbox{Tr}_q}\omega,
\cr
d(d)\;d = q^{2}\;d\;(d) +
\frac{1 - q^2}{2}\;d^2\;{\mbox{Tr}_q}\omega,
\cr
db\;d = q\;d\;db + (q^2 - 1)\;b\;d(d) +
\frac{1 - q^2}{2}\;b\;d\;{\mbox{Tr}_q}\omega,
\cr
d(d)\;b = q\;b\;d(d) +
\frac{1 - q^2}{2}\;d\;b\;{\mbox{Tr}_q}\omega,
\end{array}
\ee
\be
\begin{array}{l}
da\;b = q\;b\;da +
\frac{(q^2 - 1)}{q^2 D_{q}}\;a\;b\;(qc\;db - a\;d(d)) +
\frac{q^2 - 1}{2q^2}\;a\;b\;{\mbox{Tr}_q}\omega,
\cr
da\;d = d\;da +
\lambda\;b\;dc +
\frac{q^2 - 1}{D_{q}}\;a\;d(d\;da - \frac {1}{q}\;b\;dc) -
\frac{q^2 - 1}{2}\;a\;d\;{\mbox{Tr}_q}\omega
\cr
dc\;b = b\;dc + \frac{(q^2 - 1)}{D_{q}}\;c\;b\; (d\;da -
\frac {1}{q}\;b\;dc) - \frac{(q^2 - 1)}{2}\;c\;b\;{\mbox{Tr}_q}\omega,
\cr
dc\;d = q\;d\;dc + \frac{(q^2 - 1)}{D_{q}}\;c\;d\; (d\;da -
\frac {1}{q}\;b\;dc) - \frac{(q^2 - 1)}{2}\;c\;d\;{\mbox{Tr}_q}\omega,
\cr
db\;a = q^{-1}\;a\;db +
\frac{(q^2 - 1)}{q^2 D_{q}}\;b\;a\;(q\;c\;db - a\;d(d)) +
\frac{(q^2 - 1)}{2q^2}\;b\;a\;{\mbox{Tr}_q}\omega,
\cr
db\;c = c\;db + \frac{(q^2 - 1)}{D_{q}}\;b\;c\; (d\;da -
\frac {1}{q}\;b\;dc) - \frac{(q^2 - 1)}{2}\;b\;c\;{\mbox{Tr}_q}\omega,
\cr
d(d)\;a = a\;d(d) - \lambda\;c\;db + \frac{q^2 - 1}{D_{q}}\;
d\;a\;(d\;da - \frac {1}{q}\;b\;dc) -
\frac{(q^2 - 1)}{2}\;d\;a\;{\mbox{Tr}_q}\omega
\cr
d(d)\;c = q^{-1}\;c\;d(d) + \frac{(q^2 - 1)}{D_{q}}\;d\;c\;
(d\;da - \frac {1}{q}\;b\;dc) -
\frac{(q^2 - 1)}{2}\;d\;c\;{\mbox{Tr}_q}\omega
\end{array}
\ee

{\bf III. Differential calculus on the subgroups of the
$ GL_{q}(2,C) $ and on the their real forms.}

To obtain the differential
calculus on the Borel subgroups $B_{L}^{(2)}(C)$ and $B_{U}^{(2)}(C)$ of
the lower and the upper triangular matrices it is necessary to use in the
formulas \\(20,21,23,27,28) the surjections $\pi:  GLq(2,C) \to
B_{L}^{(2)}(C)$ such that $\pi (b) = 0, \pi (\omega^2) = 0, \pi
({\stackrel{\la}{\nabla_2}}) = 0$ and $\pi:  GLq(2,C) \to B_{U}^{(2)}(C)$
such that $\pi (c) = 0, \pi(\omega^3) = 0, \pi
({\stackrel{\la}{\nabla_3}}) = 0 $.

\noindent For the unitary matrices $ g $ and
antihermitian $\omega^{k}$:  $g^+ = g^{-1},{\omega^{k}}^+ = - \omega^{k}$.
\be
g = \left(\begin{array}{cc} a & - qD_{q}c^*\\ c & D_{q}a^* \end{array}
\right ),\; \begin{array}{cc} aa^* + q^{2}cc^* = 1, & a^*a + c^*c =
1, \\ cc^* = c^*c & D_{q}D_{q}^* = 1, \end{array} \;\;q^*=q
\ee the formulas (20,21,23,27,28) define the left differential calculus on
the $Uq(2)$ group.

\noindent To obtain the left differential calculus on the
$SL_q(2,C)$ it is necessary to suppose additional constraints $D_{q} = 1$
and $dD_{q} = D_{q}{\mbox{Tr}_q}\omega = 0$. It is possible, so the quantum
determinant $D_{q}$ commutes with the 1-forms and $dD_{q} = 0$ is
satisfied by vanishing $\bar{\omega^1} = {\mbox{Tr}_q}\omega =
\alpha(q\omega^1 + \frac {1}{q} \omega^4) = 0$.  As a consequence of these
conditions we obtained three-dimensional differential calculus, which is
independent of the parameters $\alpha$ and $\beta$.  Thus, we have the
next commutation relations on the $SL_q(2,C)$ group
\be
\begin{array}{l} \omega^1a = q^{-2}a\omega^1\cr \omega^2a
= q^{-1}a\omega^2\cr \omega^3a = q^{-1}a\omega^3 \end{array} \;\;
\begin{array}{l}
\omega^1d = q^2d\omega^1\cr
\omega^2d = qd\omega^2\cr
\omega^3d = qd\omega^3
\end{array}
\begin{array}{l}
(\omega^1)^2 = (\omega^2)^2 = (\omega^3)^2 = 0\cr
\omega^1\omega^2 + q^4\omega^2\omega^1 = 0\cr
\omega^1\omega^3 + q^{-4}\omega^3\omega^1 = 0\cr
\omega^2\omega^3 + q^{-2}\omega^3\omega^2 = 0
\end{array}
\ee
and other relations are made by interchanging $a\leftrightarrow c$,
$d\leftrightarrow b$.
The algebra of vector fields in this case has the
following form
\be
\begin{array}{lcl}
q^2{\stackrel{\la}{\nabla_1}}{\stackrel{\la}{\nabla_3}} -
q^{-2}{\stackrel{\la}
{\nabla_3}}{\stackrel{\la}{\nabla_1}} & = &
(1+q^2){\stackrel{\la}{\nabla_3}}\cr
q^2{\stackrel{\la}{\nabla_2}}{\stackrel{\la}{\nabla_1}}
- q^{-2}{\stackrel{\la} {\nabla_1}}{\stackrel{\la}{\nabla_2}}
& = &
(1 + q^2){\stackrel{\la}{\nabla_2}}\cr
{\stackrel{\la}{\nabla_3}}{\stackrel{\la}{\nabla_2}} -
q^2{\stackrel{\la}
{\nabla_2}}{\stackrel{\la}{\nabla_3}} & = &
{\stackrel{\la}{\nabla_1}}
\end{array}
\ee
\noindent The quantum group $Sp_{q}(2,C)$ is defined by condition $T^t
J_{2} = J_{2} T^{-1}$, where
\be
J_{2} = \left(\begin{array}{cc} 0 & q^{-1\over
2}\\ -q^{1\over 2} & 0 \end{array} \right )
\ee
\noindent and coincides with $SL_{q}(2,C)$ group.

 It is possible to introduce the two solutions of the differential
calculus of Wess - Zumino on the quantum plane $C_{q}(2|0)$ as the Hopf
algebra surjections $\pi:SL_q(2,C)\to B_{L}(C) = SL_{q}(2,C)\cap
B_{L}^{(2)}(C)$ such that $\pi(b) = 0$, $\pi(\omega^2) = 0$,
$\pi({\stackrel{\la}{\nabla_2}}) = 0$ and $\pi:SL_{q}(2,C)\to
B_{U}(C) = SL_{q}(2,C)\cap B_{U}^{(2)}(C)$ such that $\pi(c) = 0$,
$\pi(\omega^3) = 0$, $\pi({\stackrel{\la}{\nabla_3}}) = 0$.
\be
g_- = \left(\begin{array}{cc} a & 0 \\  c& d
\end{array}\right ) = \left(\begin{array}{cc} x & 0 \\
y&x^{-1}\end{array}\right ),g_+ = \left(\begin{array}{cc}
a & b \\
0 & d
\end{array}\right ) = \left(\begin{array}{cc}
x^{-1} & y \\
0 & x
\end{array}\right )
\ee
\noindent The quantum group $SU_{q}(2)$ = $SL_{q}(2,C)\cap U_{q}(2)$ and
the differential calculus on this group can be obtain from the
differential calculus of the $SL_{q}(2,C)$ group by the surjection
$\pi:SL_{q}(2,C)\to SU_{q}(2)$ such that $\pi(d) = a^{*}$, $\pi(b)
= - qc^*$.

There are  subgroup $U_{q}(1)\otimes U_{q}(1)$: the
diagonal subgroup of the $U_{q}(2)$ which obtained by the surjection
$\pi:U_{q}(2)\to T_{q}(2,C)$ such that $\pi(c) = 0$ .
\be
T_{q}(2,C) = \left(\begin{array}{cc}
a & 0\\
0 & D_{q}a^*
\end{array}\right )
\ee
and the diagonal subgroup of $GL_{q}(2,C)$ ($\pi:GL_{q}(2,C)\to
T_{q}(2,C)$ such that $\pi(b) = \pi(c) = 0$).
The differential calculus on this groups is the deformed calculus.

The $SU_{q}(2)$ group has the $U_{q}(1)$ subgroups: the diagonal
subgroup of $U_{q}(2)$
\be
g = \left(\begin{array}{cc}
a & 0\\
0 & a^*
\end{array}\right )
\ee
and diagonal subgroup of $SL_{q}(2,C)$
\be
g = \left(\begin{array}{cc}
a & 0\\
0 & a^{-1}
\end{array}\right )
\ee

There are, also, two subgroups $Z_{-}^{(2)}(C)$ and $Z_{+}^{(2)}(C)$
which can be considered as cosets
$Z_{-}^{(2)}(C) = GL_{q}(2,C)/B_{U}^{(2)}(C)$ and
$Z_{+}^{(2)}(C) = GL_{q}(2,C)/B_{L}^{(2)}(C)$ . The matrix elements $a, b,
d $ commute with $c$ for the $Z_{-}^{(2)}(C)$ and matrix elements $a, c,
d$ with $b$ for the $Z_{+}^{(2)}(C)$ .

The quantum subgroup
$Sp_{q}(2) = SU_{q}(2)\cap Sp_{q}(2,C)$ coincides with the $SU_{q}(2)$
group.

The quantum subgroups $O_{q}(2,C), SO_{q}(2,C)$ and
$SO_{q}(1,1)$ are the discrete group $Z_{4}$, $Z_{2}$.

The real forms of the quantum group $GL_{q}(2,C)$ and of its
subgroups and the differential calculus on this group and subgroups can be
obtained taking into account the conditions $a^+ = a$, $b^+ = b$, $c^+ =
c$, $d^+ = d$. Thus we can obtain the differential calculi on the quantum
group $GL_{q}(2,R)$ and on its quantum subgroups $B_{L}^{(2)}(R)$,
$B_{U}^{(2)}(R)$, $SL_{q}(2,R)$, $SP_{q}(2,R)$, $B_{L}(R)$,
$B_{U}(R)$ and on the real quantum plane $R_{q}(2|0)$ .

\noindent{\bf IV. WZNW model $SL_{q}(2,R)$ and $\sigma$-models on the
quantum group manifolds $SL_{q}(2,R)/U_{h}(1), C_{q}(2|9)$ and
infinitesimal transformations on the quantum group space.}

\noindent{\tt 1. Differential calculus on the $SL_{q}(2,R)$}

The matrix quantum group [4] $G = SL_{q}(2,R)$ is defined by the
q-commutation relations (C.R.) of its group parameters. Let
\be
g =
\left(\begin{array}{cc} a & b\\ c & d \end{array} \right ),\;
\begin{array}{ccc} a b = qb a , & c d = qd c, & c b =
b c \\ a c = qc a, & c d = qd c, & a d = d a + ( q -
q^{-1}) b c \end{array} \ee ${a, b, c, d}$ - hermitian, $|q|=1,
Det_{q}g = a d - qb c = 1$.\\  For any elements $g$, $g'$ $\in$
$SL_{q}(2,R)$ element $g'' = {g'}{g}$ will belong to $SL{q}(2,R)$ if $a^{k
'}a^{l} = a^{l}a^{k '}$. In the Gauss decomposition [7] \be
g= \left( \begin{array}{cc} 1&\varphi_{-}\\ 0 & 1 \end{array} \right )
 \left( \begin{array}{cc} 1 & 0\\ \varphi_{+}&1 \end{array} \right ) \left(
\begin{array}{cc}
\rho & 0\\
0 & \rho^{-1}
\end{array}
\right )=
\left(
\begin{array}{cc}
{\rho + \varphi_{-}\varphi_{+}\rho}&\varphi_{-}\rho^{-1}\\
\varphi_{+}\rho&\rho^{-1}
\end{array}
\right )
\ee
the C.R. are:
\be
\rho\varphi_{\pm} = q\varphi_{\pm}\rho,\;\;\varphi_{-}\varphi_{+} =
q^{2}\varphi_{+}\varphi_{-}
\ee
Let the quantum group is a manifolds of any possible transformations
$g' = g g_{0}$. There are two kinds of the variation: the variation
in the neighborhood of the arbitrary point of the group space $g'= g
+ dg$ and variation in the neighborhood of the unit of the group $g = 1 +
\delta g$. First variation defines the group invariants: element of the
distance between two neighboring points, element of the volume around the
point. Second variation defines the group symmetry of this invariants.
 The C.R. between the variation $dg$ and $g$ define the type of the
differential calculus.

Let $\omega = g^{-1}dg$ is the
left differential Maurer-Cartan 1-form
\be
\omega = \left(
\begin{array}{cc} \omega^{1} & \omega^{2}\\
\omega^{3} & \omega^{4}
\end{array} \right ),\;\; {\mbox{Tr}_q}\omega = q^{2}\omega^{1} +
\omega^{4} = 0
\ee
 The second order differential calculus on the $SL_{q}(2,R)$ group is
defined by the C.R.
\be
\begin{array}{ccc} \omega^{1}\rho = {1\over
q^{2}}\rho\omega^{1},& \omega^{2}\rho = {1\over q}\rho\omega^{2}, &
\omega^{3}\rho = {1\over q}\rho\omega^{3}\\
\omega^{1}\varphi_{\pm} = \varphi_{\pm}\omega^{1}, &
\omega^{2}\varphi_{\pm} = \varphi_{\pm}\omega^{2}, &
\omega^{3}\varphi_{\pm} = \varphi_{\pm}\omega^{3}
\end{array}
\ee
 The C.R. between the group parameters and their differentials are more
complicated:
\be
\begin{array}{cc}
{d\rho \rho = {1\over q^{2}}\rho d\rho}, &
{d\varphi_{+} \varphi_{+} = {1\over q^{2}}\varphi_{+} d\varphi_{+} +
(q^{4} - 1)\varphi_{+}^{3} d\varphi_{-}}\\
{d\varphi_{-} \varphi_{-} = q^{2}\varphi_{-} d\varphi_{-}}, &
{d\rho \varphi_{-} = q\varphi_{-} d\rho}
\end{array}
\ee
\be
\begin{array}{ccc}
d\varphi_{-}\varphi_{+} = q^{2}\varphi_{+}d\varphi_{-}, &
d\varphi_{+}\varphi_{-} = {1\over q^{2}}\varphi_{-}d\varphi_{+}, &
d\rho d\varphi_{-} = - qd\varphi_{-}d\rho\\
d\varphi_{-}d\varphi_{+} = - q^{2}d\varphi_{+}d\varphi_{-}, &
d\varphi_{-}\rho = {1\over q}\rho d\varphi_{-}
\end{array}
\ee
\be
\begin{array}{cc}
d\varphi_{+}\rho = {1\over q}\rho d\varphi_{+} -
q(q^{2}-1)\varphi_{+}^{2}\rho d\varphi_{-},\;\;\;
d\rho\varphi_{+} = q\varphi_{+}d\rho - q^{2}(q^{2} - 1)\varphi_{+}^{2}
\rho d\varphi_{-}\\
d\rho d\varphi_{+} + qd\varphi_{+}d\rho + q^{3}(q^{2} - 1)
\varphi_{+}^{2}d\varphi_{-}d\rho - {(q^{4} - 1)\over q^{3}}
\varphi_{+}\rho d\varphi_{-}d\varphi_{+} = 0
\end{array}
\ee
 The Maurer-Cartan 1-forms are:
\be
{\omega^{1} = \rho^{-1}d\rho+\varphi_{+}d\varphi_{-}},\;\;
{\omega^{3} = {1\over q}\rho^{2}d\varphi_{+} - q^{5}\varphi_{+}^{2}
\rho^{2}d\varphi_{-}},\;\;
{\omega^{2} = q\rho^{-2}d\varphi_{-}}
\ee
\be
\begin{array}{cc}
(\omega^1)^2 = (\omega^2)^2 = (\omega^3)^2 = 0 &\;\;\;
\omega^4 = - q^2\omega^1\\ \omega^1\omega^2 + q^4\omega^2\omega^1 =
0&\;\;\; \omega^1\omega^3 + q^{-4}\omega^3\omega^1 = 0\\ \omega^2\omega^3
+ q^{-2}\omega^3\omega^2 = 0 \end{array} \ee The left vector
fields
${\stackrel{\la}{\nabla_{k}}}$ can be obtained from the applying
the left differential to an arbitrary function on the quantum group $df =
(f{\stackrel{\la}{\partial}\over\partial a^k})da^k =
(f{\stackrel{\la}{\nabla}_{k}})\omega^{k}$.  \be
\begin{array}{cc}
{\stackrel{\la}{\nabla}} = \left(\begin{array}{cc}
{\stackrel{\la}{\nabla_{1}}}&{\stackrel{\la}{\nabla}_{2}}\\
{\stackrel{\la}{\nabla}_{3}}&{\stackrel{\la}{\nabla}_{4}}
\end{array}
\right )&\;\;\hat{\stackrel{\la}{\nabla_{1}}} =
{\stackrel{\la}{\nabla}_{1}}-q^2{\stackrel{\la}{\nabla}_{4}}
\;\;\;\;\;\;
\hat{\stackrel{\la}{\nabla_{4}}} =
{\stackrel{\la}{\nabla}_{1}}+{\stackrel{\la}{\nabla}_{4}}
\end{array} \ee
The C.R.  for vector fields have the following form \be
\begin{array}{ccc}
\rho\hat{\stackrel{\la}{\nabla}_{1}} = {1\over
q^2}\hat{\stackrel{\la}{\nabla}_{1}}\rho + \rho, &\;\;
\varphi_{-}\hat{\stackrel{\la}{\nabla}_{1}} =
\hat{\stackrel{\la}{\nabla}_{1}}\varphi_{-}, &\;\;
\varphi_{+}\hat{\stackrel{\la}{\nabla}_{1}} =
\hat{\stackrel{\la}{\nabla}_{1}}\varphi_{+}, \\ \rho{\stackrel
{\la}{\nabla}_{2}} =
{1\over q}{\stackrel{\la}{\nabla}_{2}}\rho - \varphi_{+}{\rho}^3, &
\;\;
\varphi_{-}{\stackrel{\la}{\nabla}_{2}} = {\stackrel{\la}
{\nabla}_{2}}\varphi_{-} +
{1\over q}\rho^2, &\;\; \varphi_{+}{\stackrel{\la}{\nabla}_{2}} =
{\stackrel{\la}{\nabla}_{2}}\varphi_{+} + q\varphi_{+}^2\rho^2\\
\rho{\stackrel{\la}{\nabla}_{3}} = {1\over
q}{\stackrel{\la}{\nabla}_{3}}\rho, &\;\;
\varphi_{-}{\stackrel{\la}{\nabla}_{3}} =
{\stackrel{\la}{\nabla}_{3}}\varphi_{-}, &\;\;
\varphi_{+}{\stackrel{\la}{\nabla}_{3}} =
{\stackrel{\la}{\nabla}_{3}}\varphi_{+} + q\rho^{-2}, \\ \end{array}
\ee
\be \begin{array}{c}
q^2\hat{{\stackrel{\la}{\nabla}_{1}}}{\stackrel{\la}
{\nabla}_{3}} - {1\over q^2}{\stackrel{\la}{\nabla}_{3}}
\hat{{\stackrel{\la}{\nabla}_{1}}} =
(q^2 + 1){\stackrel{\la}{\nabla}_{3}},\cr
q^2{\stackrel{\la}{\nabla}_{2}}\hat{{\stackrel{\la}
{\nabla}_{1}}} - {1\over q^2}\hat{\stackrel{\la}{\nabla}_{1}}
{\stackrel{\la}{\nabla}_{2}} = (q^2 + 1){\stackrel{\la}
{\nabla_{2}}},\cr {\stackrel{\la}{\nabla}_{3}}{\stackrel{\la}
{\nabla}_{2}} - {1\over q^2}{\stackrel{\la}{\nabla}_{2}}{\stackrel
{\la}{\nabla}_{3}} = \hat{\stackrel{\la}{\nabla}_{1}}
\end{array}
\ee
and $\stackrel{\la}{\nabla}_{k}$ have the form
\be
\hat{\stackrel{\la}{\nabla}_{1}} = {{\stackrel{\la}{\partial}}
\over \partial\rho}\rho,\;\;{\stackrel{\la}{\nabla}_{2}} = {1\over
q}{{\stackrel{\la}{\partial}}\over \partial\varphi_{-}}\rho^2 -
{{\stackrel{\la}{\partial}}\over \partial\rho}\varphi_{+}\rho^3 +
q{{\stackrel{\la}{\partial}}\over
\partial\varphi_{+}}\varphi_{+}^2\rho^2,\;\;{\stackrel{\la}
{\nabla}_{3}} = q{{\stackrel{\la}{\partial}}\over
\partial\varphi_{+}}\rho^{-2} \ee The left
vector fields and the left derivatives act on the any function of the
group parameters from the right side.

\noindent{\tt 2.WZNW model on the
$SL_{q}(2,R)$ group.}

\noindent The existing of the quantum group structure in the
WZNW model was shown in [8,9]. The $\sigma$-models on the quantum group
spaces was considered in [7,10-12]. To construct the WZNW model with
$SL_{q}(2,R)$ group symmetry, we consider the space ${M^{(1,1)}\oplus
SL_{q}(2,R)}$, where base space $M^{(1,1)}$ is the commutative
(nondeformed) space.  The element of the volume in $M^{(1,1)}$ space,
which is the invariant of $SL_{q}(2,R)$ group, is
\be
{{\mbox{Tr}_q}[\omega(d)\wedge dz^{\mu}][\omega(d)\wedge dz^{\mu}]\over
{2\epsilon ^{\lambda\rho} dz^{\lambda}\wedge dz^{\rho}}} =
{\mbox{Tr}_q}(\omega_{\mu}\omega^{\mu}) d^{2}z ,
\ee
\noindent where
$\omega(d) = \omega_{\mu} dz^{\mu}, z^{\mu}\epsilon M^{1,1}, \mu = 1,2.$
For any $2\times 2$ matrix $A, {\mbox{Tr}_q} A = q^{2}A^{1} + A^{4}$.
As a result
we have
$$ {\mbox{Tr}_q}(\omega_{\mu}\omega^{\mu}) = q^{5}[2]_{q}
\rho^{-2}\partial_{\mu}\rho\partial^{\mu}\rho+q^{5}
[2]_{q}\rho^{-1}\varphi_{+}(\partial_{\mu}\varphi_{-}
\partial^{\mu}\rho + {1\over q}
\partial_{\mu}\rho\partial^{\mu}\varphi_{-}) +
$$
$$
(\partial_{\mu}\varphi_{-}\partial^{\mu}\varphi_{+} +
q^{2}\partial_{\mu}\varphi_{+}\partial^{\mu}\varphi_{-}) -
q^{2}(q^{4} - 1)\varphi_{+}^{2}\partial_{\mu}\varphi_{-}
\partial^{\mu}\varphi_{-}
$$
\noindent\underline{ The C.R. are now in the same space-time point} and
$d\rho = \partial_{\mu}\rho dz^{\mu}, d\varphi_{\pm} =
\partial_{\mu}\varphi_{\pm} dz^{\mu}$ and
$[n]_{q} = {{q^{n} - q^{-n}}\over{q - q^{-1}}}$. The Wess-Zumino term
\be
Tr_{q}(\omega (d)\wedge\omega (d)\wedge \omega (d)) =
{q[2]_{q}[3]_{q}\over 6}\varepsilon^{\mu\nu\lambda}
\partial_{\lambda}(\rho^{-1}\partial_{\mu}\rho\partial_{\nu}\varphi_{-}
\varphi_{+}) d^{3}z
\ee
\noindent is the total derivative.
Finally, the WZNW-action  with the $SL_{q}(2,R)$ quantum group symmetry
describes the 2-dimensional relativistic string in the background gravity
and antisymmetric fields
\be S[\rho,\varphi_{-},\varphi_{+}] = {k\over
4\pi}\int d^{2}z(G_{AB}\partial_{\mu}X^{A}\partial^{\mu}X^{B} +
B_{AB}\varepsilon_{\mu\nu}\partial^{\mu}X^{A}\partial^{\nu}X^{B}) ,
\ee
where $X^{A} = (\rho,\varphi_{-},\varphi_{+})$ and the background gravity
and antisymmetric fields have the following form:  \be
G_{AB} =
\left(
\begin{array}{ccc}
q^{5}[2]_{q}\rho^{-2}&q^{4}[2]_{q}\rho^{-1}\varphi_{+} & 0\\
q^{5}[2]_{q}\rho^{-1}\varphi_{+} & -q^{2}(q^{4} - 1)\varphi_{+}^{2} & 1\\
0 & q^{2} & 0
\end{array}
\right )\ee
\be B_{AB} =
{q^{3}[2]_{q}[3]_{q}\over 6}\varphi_{+}\rho^{-1}
\left(
\begin{array}{ccc}
0 & 1 & 0\\
-1 & 0 & 0\\
0 & 0 & 0
\end{array}
\right )
\ee
 The group symmetry of this model is $SL_{q}(2,R)\otimes SL_{q}(2,R)$,
because under the left multiplication on the group
$g' = g_{0}g$ the differential forms of Cartan are invariant, $\omega
^\prime = \omega$, and under the right multiplication $g' = gg_{0}$
the differential forms are covariant, $\omega ^{\prime} =
{g_{0}^{-1}}\omega g_{0}$. But $Tr_{q}A$ is invariant of the
transformation $A' = g_{0}^{-1}Ag_{0}$, because the elements of matrix $A$
commute with the elements of matrix $g_{0}$, by definition of the quantum
group. Therefore, this model describes the spontaneous breaking of the
$SL_{q}(2,R)\otimes SL_{q}(2,R)$ symmetry to the $SL_{q} (2,R)$ one.

\noindent{\tt 3.$\sigma$-model
on the $SL_{q}(2,R)/U_{h}(1)$ group.}

\noindent Let us consider the spontaneous
breaking symmetry in the $\sigma-$ model with the \\ $SL_{q}(2,R)/U_{h}(1)$
group symmetry.  Let $G = KH$, K-coset, H-subgroup. The Maurer-Cartan
1-forms
\be
k^{-1}dk = \left( \begin{array}{cc}
q^{2}\varphi_{+}d\varphi_{-}&d\varphi{-}\\
d\varphi_{+} - q^{2}\varphi_{+}^{2}d\varphi_{-} &
- \varphi_{+}d\varphi_{-}
\end{array}
\right ) = \omega + \theta ,
\ee
where $\omega\epsilon K, \;\;\theta\epsilon H$ and the coset elements
$\varphi_{\pm}$ commute with the subgroup parameter $\rho$ and satisfy to
C.R. of $SL_{q}(2,R)$ group among themselves.  There is a question: how
does coset and subgroup separate from $k^{-1}dk$ ?  In opposite to the
classical case, there is the 3-parametric family of the $U(1)$
subgroups. The Lagrangian has
the following form:
\be
L_{n} = {1\over 2}Tr_{q}(\omega_{\mu}\omega^{\mu}) = {(q^{4}+1)\over
4q^{4}}(\partial_{\mu}\varphi_{-}\partial^{\mu}
\varphi_{+}+q^{2}\partial_{\mu}\varphi_{+}\partial^{\mu}\varphi_{-})
- c_{n}(q)\varphi_{+}^{2}\partial_{\mu}\varphi_{-}\partial
^{\mu}\varphi_{-},
\ee
where $c_{n}(q)$ depends on the choice of a subgroup. There are three
most interesting examples.

\noindent{\tt 1)}. Nondeformed $U(1)$ subgroup : $c_{1} = {{2q^{4} -
q^{2}+1}\over 2}$ \be \omega = \left( \begin{array}{cc}
(q^{2} -
 1)\varphi_{+}d\varphi_{-}&d\varphi_{-}\\ d\varphi_{+} -
q^{2}\varphi_{+}^{2}d\varphi_{-} & 0 \end{array} \right ),\;\;\;\theta =
\varphi_{+}d\varphi_{-} \left( \begin{array}{cc} 1 & 0\\ 0 & - 1
\end{array}
\right )
\ee
\noindent The algebra symmetry of this Lagrangian is defined by the
Maurer-Cartan equations:
\be d\theta = - \left( \begin{array}{cc}
q^{-2} & 0\\
0 & 1
\end{array}
\right )
\omega\omega + (q^{2}-1)\omega\theta , d\omega = - \left(
\begin{array}{cc} {q^{2} - 1}\over{q^{2}} & 0\\ 0 & 0 \end{array} \right )
\omega\omega - q^{3}[2]_{q}\omega\theta , \theta\omega = q^{4}\omega\theta
\ee
\noindent The C.R. between
the coset and the subgroup forms are common for all of the examples
\be
\begin{array}{cc} {\omega ^{1}\omega ^{3} + q^{4}\omega ^{3}\omega ^{1}
= 0}, & {\omega ^{2}\omega ^{3} + q^{2}\omega ^{3}\omega ^{2} = 0}\cr
{\omega ^{1}\omega ^{3} + q^{4}\omega ^{3}\omega ^{1} = 0}, & {\omega
^{4}\omega ^{3} + q^{4}\omega ^{3}\omega ^{4} = 0} \end{array}
\ee

\noindent{\tt 2)}. Classical coset structure: $c_{2} = {{q^{6} + 1}\over
4}$ \be \omega = \left( \begin{array}{cc} 0 & d\varphi_{-}\\
d\varphi_{+} -
q^{2}\varphi_{+}^{2}d\varphi_{-} & 0 \end{array} \right ),\;\;\;\theta =
\varphi_{+}d\varphi_{-} \left(\begin{array}{cc}
q^{2} & 0\\
0 & - 1
\end{array}
\right )
\ee
\be
d\theta = - \omega\omega,\;\;
d\omega = - \omega\theta - \theta\omega ,\;\theta\omega =
q^{2}\omega\theta
\ee

\noindent{\tt 3)}. There is one of the examples of the 2-
parametric family $U_{q}(1)$ subgroups:
\be
\begin{array}{ccc}c_{3} = {{2q^{4} - q^{2} + 1}\over 2q^{2}}, &
\omega = \left( \begin{array}{cc} {(q^{2} - 1)\over
q^{2}}\varphi_{+}d\varphi_{-} & d\varphi_{-}\\
d\varphi_{+}-q^{2}\varphi_{+}^{2}d\varphi_{-} & 0
\end{array}
\right ),&\theta = {1\over q^{2}}\varphi_{+}d\varphi_{-}
\left(\begin{array}{cc}
1 & 0\\
0 & - q^{2}
\end{array}
\right )
\end{array}
\ee
\be
d\theta = -
\left(
\begin{array}{cc}
q^{-4} & 0\\
0 & 1
\end{array}
\right )
\omega\omega + (q^{4}-1)\omega\theta,\;
d\omega = -
\left(
\begin{array}{cc}
{{q^{4} - 1}\over{q^{4}}} & 0\\
0 & 0
\end{array}
\right )
\omega\omega - q^{4}[2]_{q}\omega\theta,\theta\omega = q^{6}\omega\theta
\ee
Why we have obtained different algebras of a symmetry for the same
subgroup?
That is possible because we can use the different map from the algebra to
the group, for example:
\be g = \exp(\varphi_{-}\tau_{+})\exp(\varphi_{+}\tau_{-})
\exp({\ln\rho}\tau_{3}) ,
\ee
where $\tau$ are the Pauli matrices -- the fundamental representation of
the $U_{q}(SL(2,R)$ algebra.The group stability of the vacuum is $U(1)$.
In the another parametrization
\be g = \exp(\varphi_{-}\tau_{+})\exp(\varphi_{+}\tau_{-})
(1 - {(q^{2} - 1)\over
q^{2}}\nabla_{3})^{\ln\rho\over \ln q^{-2}},\; \nabla_{3} = \left(
\begin{array}{cc}
1 & 0\\
0 & - q^{2}
\end{array}
\right )
\ee
the group stability of the vacuum is $U_{q}(1)$.
 The group symmetry of this Lagrangians is $SL_{q}(2,R)$ spontaneously
broken to $U_{p}(1),\;\;p = q^{\pm 2n}, n = 0,1... $.  Under the left
multiplication on the group $g' = g_{0}g$, the differential form
$g^{-1}dg = h^{-1}(\omega + \theta )h = {g'} ^{-1}dg'$.  Therefore,
$\omega^{\prime} + \theta^{\prime} = h'h^{-1}(\omega + \theta )h{h'}^{-1}$
.  Again, the decomposition on the coset and the subgroup forms is not
unique after the transformation . The group transformation can transform
the Lagrangian with the $U_{p_{1}}(1)$ subgroup of the vacuum stability to
the Lagrangian with the $U_{p_{2}}(1)$ subgroup.

\noindent{\tt 4.Variational
calculus on the $SL_{q}(2,R)$ group.}

\noindent It is possible to obtain  the
variational calculus on the group by two ways: from the C.R. between the
left vector fields and group parameters and from the infinitesimal
transformations on the group. Let us multiply the C.R. (48) between
$\nabla_{n}$ and group parameters on the parameters of transformation
$R^n$. The form of the infinitesimal transformations of the group
parameters is obtained under the requirement \be
[X_{A},{\stackrel{\la}{\nabla}_{n}}R^n] = X_{A}{\stackrel{\la}
{\delta}} = X_{A}\sum_{n = 1}^{3}{\stackrel{\la}{\delta}_{R_{n}}},
\;\;\;X_{A} = (\rho, \varphi_{-}, \varphi_{+}),\;\;\; [ A , B ] = AB-BA
\ee
By imposing the
C.R. between the parameters of infinitesimal transformations and group
parameters \be \begin{array}{ccc} \rho R^{1} = q^2
 R^{1}\rho, &\;\;\varphi_{-}R^{1} = R^{1}\varphi_{-}, &\;\;\varphi_{+}R^{1}
= R^{1}\varphi_{+}\\ \rho R^{2} = q R^{2}\rho, &\;\;\varphi_{-}R^{2} =
R^{2}\varphi_{-}, &\;\;\varphi_{+}R^{2} = R^{2}\varphi_{+}\\ \rho R^{3} = q
R^{3}\rho, &\;\;\varphi_{-}R^{3} = R^{3}\varphi_{-}, &\;\;\varphi_{+}R^{3}
= R^{3}\varphi_{+} \end{array} \ee we obtain the infinitesimal
transformation of the group parameters
\be
\begin{array}{ccc}\rho{\stackrel{\la}{\delta}} = \rho R^{1} -
\varphi_{+}\rho^3R^{2}, &\;\;\varphi_{-}{\stackrel{\la}{\delta}} =
{1\over q}\rho^2 R^{2}, &\;\;\varphi_{+}{\stackrel{\la}{\delta}} =
q\varphi_{+}^2\rho^2 R^{2} + q\rho^{-2}R^{3} \end{array} \ee  In
the terms of the components ${\stackrel{\la}{\delta}_{R_{n}}}$ the
infinitesimal transformations have the following form:  \be
 \begin{array}{ccc} \rho{\stackrel{\la}{\delta}_{R^{1}}} = \rho
 R^{1}, & \rho{\stackrel{\la}{\delta}_{R^{2}}} = -
 \varphi_{+}{\rho}^3 R^{2}, & \rho{\stackrel{\la}{\delta}_{R^{3}}}
 = 0\\ \varphi_{-}{\stackrel{\la}{\delta}_{R^{1}}} = 0, &
 \varphi_{-}{\stackrel{\la}{\delta}_{R^{2}}} = {1\over q}{\rho}^2
 R^{2}, & \varphi_{-}{\stackrel{\la}{\delta}_{R^{3}}} = 0\\
\varphi_{+}{\stackrel{\la}{\delta}_{R^{1}}} = 0, &
\varphi_{+}{\stackrel{\la}{\delta}_{R^{2}}} = q{\varphi_{+}}^2
{\rho}^2 R^{2}, & \varphi_{+}{\stackrel{\la}{\delta}_{R^{3}}} =
q\rho^{-2} R^{3}\end{array} \ee
\noindent We postulate the $U_{q}SL(2,R)$
algebra of the vector fields in the form, which it is common for a boson
theory and a supersymmetric theory \be \begin{array}{c} [{\stackrel{\la}
{\nabla} _{1}}R^{1},{\stackrel{\la}{\nabla} _ {2}}R^{2}] =
(q^2+1){\stackrel{\la}{\nabla} _{2}}R^{1}R^{2},\cr
[{\stackrel{\la}{\nabla}_{3}}R^{3},{\stackrel{\la}
{\nabla} _{1}}R^{1}] = (q^2 + 1){\stackrel{\la}{\nabla}
_{3}}R^{3}R^{1},\cr [{\stackrel{\la}{\nabla} _{2}}R^{2},{\stackrel
{\la}{\nabla}_{3}}R^{3}] = {\stackrel{\la}{\nabla} _{1}}
R^{2}R^{3} \end{array}
\ee
It is possible, if we suppose the following C.R. between the infinitesimal
parameters and the vector fields:
\begin{equation}
\begin{array}{ccc}
R^{1}{\stackrel{\la}{\nabla}_{2}} = \frac
{1}{q^4}{\stackrel{\la}{\nabla}_{2}}R^{1} &
R^{2}{\stackrel{\la}{\nabla}_{1}} = q^{4}{\stackrel{\la}{\nabla}_{1}}R^{2}
+ q^{2}(q^{4} - 1)R^{2} & R^{1}{\stackrel{\la}{\nabla}_{3}} =
q^{4}{\stackrel{\la}{\nabla}_{3}}R^{1}\\
R^{3}{\stackrel{\la}{\nabla}_{1}} = \frac
{1}{q^4}{\stackrel{\la}{\nabla}_{1}}R^{3} + \frac {q^{4} - 1}{q^4}R^{3}&
R^{3}{\stackrel{\la}{\nabla}_{2}} = {\stackrel{\la}{\nabla}_{2}}R^{3}&
R^{2}{\stackrel{\la}{\nabla}_{3}} = {\stackrel{\la}{\nabla}_{3}}R^2
\end{array}\end{equation}
It means, that we have the same C.R. between basis of 1-forms and basis of
vector fields
\begin{equation}
\begin{array}{ccc}
\omega^{1}{\stackrel{\la}{\nabla}_{2}} = \frac
{1}{q^4}{\stackrel{\la}{\nabla}_{2}}\omega^{1} &
\omega^{2}{\stackrel{\la}{\nabla}_{1}} =
q^{4}{\stackrel{\la}{\nabla}_{1}}\omega^{2} + q^{2}(q^{4} - 1)\omega^{2} &
\omega^{1}{\stackrel{\la}{\nabla}_{3}} =
q^{4}{\stackrel{\la}{\nabla}_{3}}\omega^{1}\\
\omega^{3}{\stackrel{\la}{\nabla}_{1}} = \frac
{1}{q^4}{\stackrel{\la}{\nabla}_{1}}\omega^{3} + \frac {q^{4} -
1}{q^4}\omega^{3}& \omega^{3}{\stackrel{\la}{\nabla}_{2}} =
{\stackrel{\la}{\nabla}_{2}}\omega^{3}& R^{2}{\stackrel{\la}{\nabla}_{3}}
= {\stackrel{\la}{\nabla}_{3}}\omega^2 \end{array}\end{equation}
\noindent The same result we can
obtain from the right infinitesimal multiplication on the group $g' = g
g_{0}$, where $g_{0} = 1 + \delta g_{0}$.  For \be \delta g_{0} =
\left(\begin{array}{cc} R^{1} & R^{2}\\ R^{3} & -q^2R^{1} \end{array}
\right ) \ee we see, that $dg = g\delta g_{0}$ and C.R. for $R^n$ have the
following form:  \be \begin{array}{ccc} R^1 R^2 = q^4 R^2 R^1, & R^3 R^1 =
q^4 R^1 R^3, & R^3 R^2 = q^2 R^2 R^3 \end{array} \ee \noindent
simultaneously with the condition $R^{4} = - q^2 R^{1}$.  The C.R.  of the
variational calculus \be \begin{array}{cc} {(\rho{\stackrel{\la}{\delta}})
\rho = {1\over q^{2}}\rho (\rho{\stackrel{\la}{\delta}})}, & {(\varphi_{+}
{\stackrel{\la}{\delta}})\varphi_{+} = {1\over q^{2}}\varphi_{+}
(\varphi_{+}{\stackrel{\la}{\delta}}) +
(q^{4} - 1)\varphi_{+}^{3} (\varphi_{-}{\stackrel{\la}{\delta}})}\\
{(\varphi_{-}{\stackrel{\la}{\delta}}) \varphi_{-} = q^{2}\varphi_{-}
(\varphi_{-}{\stackrel{\la}{\delta}})}, & {(\rho{\stackrel{\la}
{\delta}})\varphi_{-} = q\varphi_{-} (\rho{\stackrel{\la}{\delta}})}
\end{array}
\ee
\be
\begin{array}{c} (\varphi_{-}{\stackrel{\la}{\delta}})\rho = {1\over
q}\rho(\varphi_{-}{\stackrel{\la}{\delta}}) \end{array} \ee
\be
\begin{array}{cc} (\varphi_{+}{\stackrel{\la}{\delta}})\rho = {1\over
q} \rho (\varphi_{+}{\stackrel{\la}{\delta}}) -
q(q^{2} - 1)\varphi_{+}^{2} \rho
(\varphi_{-}{\stackrel{\la}{\delta}}), &
(\varphi_{-}{\stackrel{\la}{\delta}})\varphi_{+} = q^{2}\varphi_{+}
(\varphi_{-}{\stackrel{\la}{\delta}})\\
(\rho{\stackrel{\la}{\delta}})\varphi_{+} = q\varphi_{+}
(\rho{\stackrel
{\la}{\delta}}) - q^{2}(q^{2} - 1)\varphi_{+}^{2} \rho
(\varphi_{-}{\stackrel{\la}{\delta}}), &
(\varphi_{+}{\stackrel{\la}{\delta}})\varphi_{-} = {1\over
q^{2}}\varphi_{-}(\varphi_{+}{\stackrel{\la}{\delta}})
\end{array}
\ee
\noindent are consistent
with the C.R.  (39) and are coincide with the C.R. of the differential
calculus (42,43).  The $U_{q}(SL(2,R))$ algebra in the form (71) is the
condition of the compatibility of the relations (69,70)
\be
\begin{array}{c}
X^A[{\stackrel{\la}{\delta}_{R^1}}
,{\stackrel{\la}{\delta_{R^2}}}
] = X^{A}(q^2+1){\stackrel{\la}{\delta}_{R^{1}R^2}}\\
X^A[{\stackrel{\la}{\delta}_{R^3}}
, {\stackrel{\la}{\delta}_{R^1}}
]
= X^{A}(q^2+1){\stackrel{\la}{\delta}_{R^{3}R^1}}\\
X^A[{\stackrel{\la}{\delta}_{R^2}}
,{\stackrel{\la}{\delta_{R^3}}}
] = X^{A}{\stackrel{\la}{\delta}_{R^{2}R^3}}
\end{array} \ee

\noindent{\tt 5.Equations of motion}

\noindent We use the extremum principle of the action to obtain the
equations of motion  and
we must to commute the variations of fields and their derivatives on the
right or on the left side. We can use both variation $dX^A$ and
$X^A{\stackrel{\la}{\delta}}$ to do this. The C.R.  of the
differential
calculus on the $SL_{q}(2,R)$ group are insufficient to do ordering.
Therefore, we need in the differential calculus on the Lagrangian
manifolds $(\rho ,\varphi_{\pm},\dot\rho ,\dot{\varphi_{\pm}}
,\acute{\rho} , \acute{\varphi_{\pm}} )$.  This is not the quantum group
manifold and we can not use the formalism of 1-forms. We can require, that
the Lagrangian equation of motion be coincident with the conservation law
$\partial_{\mu}{\omega}^{\mu}=0$ for Lagrangian with $SL_{q}(2,R)$ group
symmetry. At last, we can investigate the 1- dimensional $\sigma$- models.
 The variational calculus is more suitable to obtain the equations of
motion. The C.R. between the $X^A$, $\dot{X^A}$, $\acute{X^A}$ and $R^n$,
$\dot{R^n}$, $\acute{R^n}$ can obtain by differentiating the relations
(68).
\be \begin{array}{ccc} \dot{\rho}R^1 = q^2 R^1
\dot{\rho} &\;\dot{\rho}\dot{R^1} = q^2\dot{R^1}\dot{\rho}&\;
\rho\dot{R^1} = \dot{R^1}\rho\\
\dot{\rho}R^2 =
qR^2\dot{\rho} &\;\dot{\rho}\dot{R^2} = q\dot{R^2}\dot{\rho} &\;
\rho\dot{R^2} = q\dot{R^2}\rho\\
\dot{\rho}R^3 =
qR^3\dot{\rho}&\;\dot{\rho}\dot{R^3} = q\dot{R^3}\dot{\rho} &\;
\rho\dot{R^3} = q\dot{R^3}\rho
\end{array}
\ee
 The derivatives of $\varphi_{\pm}$ commute with the derivatives of
$R^n$.

\noindent{\tt 6. One dimensional $\sigma$- model on the quantum plane
$C_{q}(2|0)$.}

\noindent The differential calculus on the $C_{q}(2|0)$ is
coincide with the differential calculus on the Borel subgroup of
$SL_{q}(2,C)$ and can be obtained from the differential calculus on the
$SL_{q}(2,C)$ by surjection:  $\pi$:$SL_{q}(2,C)\rightarrow B_{L}$ such
that $\pi(b) = 0$.  \be g = \left( \begin{array}{cc} x & 0\\
y & x^{-1}
\end{array} \right ),\;\; \begin{array}{cc} {xy = qyx}&\dot yy =
q^{-2}y\dot y \\ \dot xx = q^{-2}x\dot x & \dot xy = q^{-1}y\dot x
\end{array} \ee \be \omega = \left(
\begin{array}{cc} x^{-1}dx & 0\\ xdy-qydx &-q^{2}x^{-1}dx \end{array}
\right ),\; \dot yx = q^{-1}x\dot y - {(q^{2} - 1)\over q^{2}}y\dot x
\ee In term of the variables $\rho,\varphi_{\pm}$
\be g = \left( \begin{array}{cc} \rho & 0\\
\varphi_{+}\rho&\rho^{-1} \end{array} \right );\;\; \omega = \left(
\begin{array}{cc} \rho^{-1}d\rho&0\\ {1\over
q}\rho^{2}d\varphi_{+} & - q^{2}\rho^{-1}d\rho \end{array} \right )
\ee
\be
L = {1\over 2}{\mbox{Tr}_q}(\omega_{\mu}\omega^{\mu}) = {q^{4}(q^{2} + 1)
\over 2}
\rho^{-2}{\dot\rho}^{2}
\ee
The equation of motion  ${\dot\omega}^{1} =
 \rho^{-1}\ddot{\rho} - q^{2}\rho^{-1}{\dot\rho}^{2} = 0$ will
coincide with Lagrangian equation , if we impose the C.R.
$d\dot\rho\dot\rho = {1\over
q^{2}}\dot\rho d\dot\rho$.
The classical solution of this equation is
\be \rho = \alpha\exp(\beta t),\;\;\
\alpha\beta = q^{2}\beta\alpha
\ee
and C.R. in the different time points are
\be \rho(t)\rho(t^{\prime}) = \rho(q^{2}t^{\prime})
\rho({1\over q^{2}}t),\;\;
\rho(t)\rho(t^{\prime}) = \exp[q^{2}(q^{2} - 1)\beta
	[(t - t^{\prime})]\rho(t^{\prime})\rho(t)
\ee
\noindent There are
$4\times 4$ matrix representations of $\alpha, \beta$ such,that
$det_{q}\alpha = 0$ or $det_{q}\beta = 0$.Therefore, we can
rewrite this Lagrangian as a $4\times 4$ matrix model for the commuting
fields.
In conclusion, we note that 2-dimensional $\sigma$-model on
the guantum plane
\be
 L = {q^{4}(q^{2} + 1)\over
2}\rho^{-2}\partial_{\mu}\rho\partial^{\mu}\rho
\ee
leads to the C.R.
$d\acute\rho\acute\rho = {1\over
q^{2}}\acute\rho d\acute\rho$ and the equation of motion
$\partial_{\mu}\partial^{\mu}\rho - q^{2}\rho^{-1}\partial_{\mu}
\rho\partial^{\mu}\rho = 0,\;\;\mu = 1,2$.

I would like to thank M. Arik and R. Mirkasimov for kind
hospitality during the conference "Quantum groups, deformations and
contractions" and V.  Dobrev, J.  Lukierski, F. Muller-Hoissen for
stimulating discussions.

This work was supported in part by the Ukrainian Ministry of Science and
technology, grant N 2.5.1/54, by grants INTAS 93-633 (Ext) and INTAS
93-127 (Ext).

\begin{center}
{\bf References}
\end{center}
\begin{enumerate}
\item V. P. Akulov, V. D. Gershun, A. I. Gumenchuk, {\it JETP Lett.},
{\bf 58} (1993), 474
\item V. Akulov, V. Gershun, {\it q-alg 9509030}
\item V. P. Akulov, V. D. Gershun,{\it Proc. Int.  conf. on math,
phys.  (Rahov, Ukraine, 1995)}, (in press)
\item L. Faddeev, N. Reshetikhin, L. Takhtajan, {\it
Algebra i Analiz} { \bf 1} (1989), 178
\item A.Schirrmacher, J.Wess, B.Zumino, {\sl Z.Phys.} {\bf C49} (1990), 317
\item S. L. Woronowicz, {\it Comm.  Math.  Phys.} {\bf 122} (1989),
125
\item V. P. Akulov, V. D. Gershun, A. I. Gumenchuk, {\it JETP Lett.},
{\bf 56} (1992), 177
\item L. D. Fadeev, {\it Commun. Math, Phys.}, {\bf 132} (1990), 131
\item A. Alekseev and S. Shatashvili, {\it Commun. Math.
Phys.}, {\bf 133} (1990), 353
\item I. Y. Aref'eva and I. V. Volovich,
{\it Phys. Lett.}, {\bf 264B} (1991), 62
\item Y. Frishman, J. Lukierski, W. J. Zakrzewski, {\it
J.Phys.A:Math.Gen}, {\bf 26} (1993), 301
\item
V. D. Gershun, {\it Proc.  Xth Int. conf. "Problems of Quantum Field
Theory" (Alushta, Ukraine,1996)}, JINR, Dubna, (1996), 119
\item V. D. Gershun {\it q-alg 9709022}
\end{enumerate}

\end{document}